\documentclass[conference]{IEEEtran}
\ifCLASSINFOpdf
   \usepackage[pdftex]{graphicx}
   \DeclareGraphicsExtensions{.pdf,.jpeg,.png}
\else
\fi
\usepackage[cmex10]{amsmath}
\begin{document}
\title{Applications of AI for Magic Squares}

\author{\IEEEauthorblockN{Jared Weed}
\IEEEauthorblockA{Department of Mathematical Sciences\\
Worcester Polytechnic Institute\\
Worcester, Massachusetts 01609-2280\\
Email: jmweed@wpi.edu}}


%Telephone: (800) 555--1212, Fax: (888) 555--1212}


\maketitle

\begin{abstract}
%\boldmath
In recreational mathematics, a \emph{normal} magic square is an $n \times n$ square matrix whose entries are distinctly the integers $1 \ldots n^2$, such that each row, column, and major and minor traces sum to one constant $\mu$. It has been proven that there are 7,040 fourth order normal magic squares and 2,202,441,792 fifth order normal magic squares [4], with higher orders unconfirmed [3]. Previous work related to fourth order normal squares has shown that symmetries such as the dihedral group exist [5] and that (under certain conditions)  normal magic squares can be categorized into four distinct subsets. [2][6] With the implementation of an efficient backtracking algorithm along with supervised machine learning techniques for classification, it will be shown that the entire set of fourth order normal magic squares can be generated by expanding the symmetry groups of 95 asymmetric parents. Discussion will suggest that methods employed in this project could similarly apply to higher orders.
\end{abstract}


% \fi
\IEEEpeerreviewmaketitle

\section{Introduction}Magic squares, being a recreation mathematical topic in nature, have been the subject of entertainment and interest to many mathematicians and math-enthusiasts alike for hundreds of years. A \emph{normal} magic square is an $n \times n$ square matrix whose entries are distinctly the integers $1 \ldots n^2$, such that each row, column, and major and minor traces sum to one constant $\mu$. Magic squares of order 3 and 4 (and even higher orders) have appeared in paintings, literature, and artifacts dated as far back as 650 BC. One of the first magic squares found in ancient Chinese literature is the \textbf{Lo Shu square}, which is told to have been painted on the shell of a sea turtle:

$$\begin{array}{|c|c|c|}
\hline 4 & 9 & 2 \\
\hline 3 & 5 & 7 \\
\hline 8 & 1 & 6 \\
\hline
\end{array}$$

One of the most notable magic squares in western civilization is the order 4 magic square in Albrecht Durer's engraving \emph{Melencolia I}:
$$\begin{array}{|c|c|c|c|}
\hline 16 & 3 & 2 & 13 \\
\hline 5 & 10 & 11 & 8 \\
\hline 9 & 6 & 7 & 12 \\
\hline 4 & 15 & 14 & 1 \\
\hline
\end{array}$$
The Durer magic square has many fascinating properties. Aside the rows and columns (and major/minor diagonal) summing to the magic number 34, the sum can also be found in each 2 by 2 quadrant,
\includegraphics[width=1.0\linewidth]{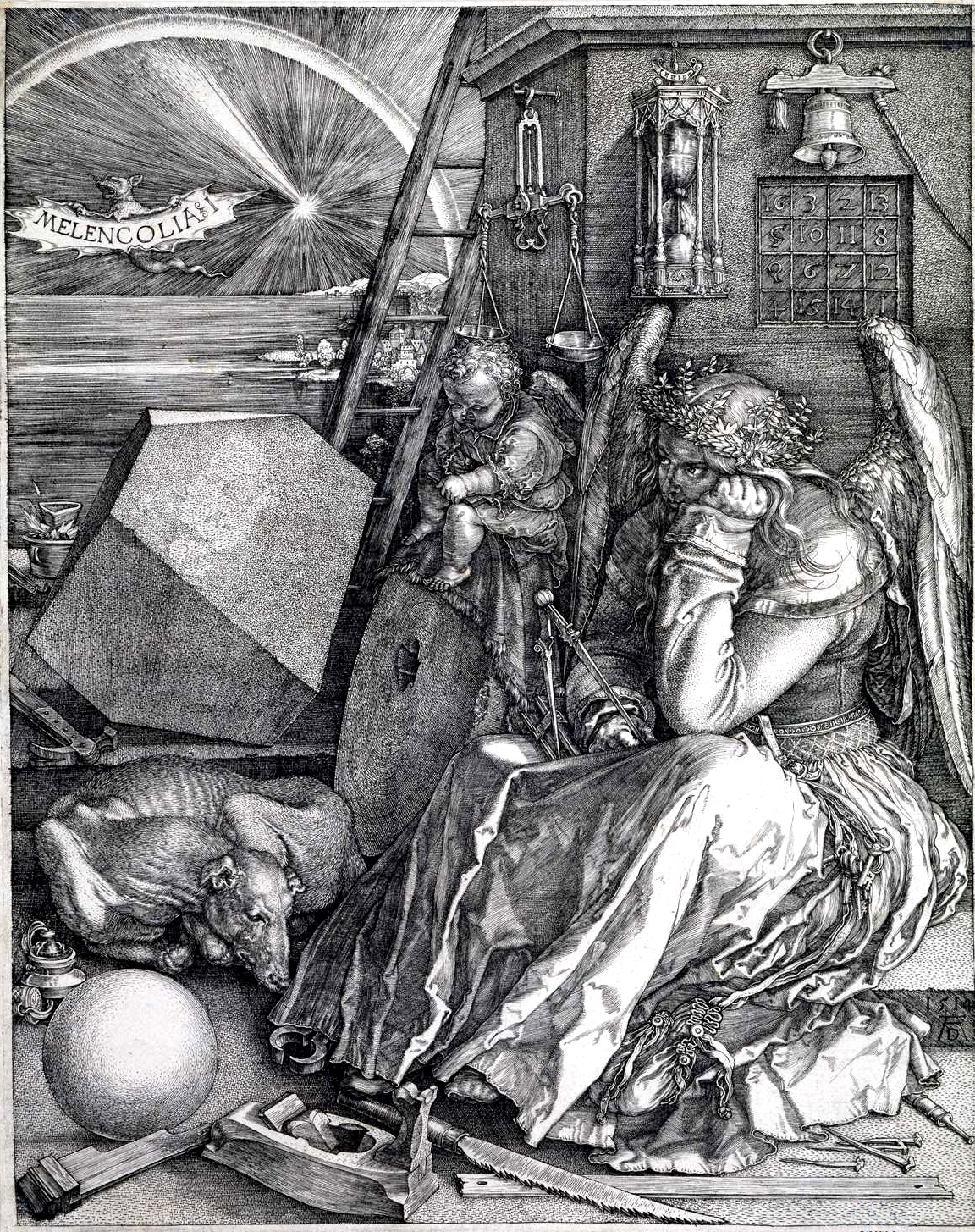} the center 2 by 2 square, the four outer corner squares, and the four outer corner squares of each 3 by 3 subsquare. Furthermore, the sum can be found clockwise from the corners ($3 + 8 + 14 + 5$) and counter-clockwise ($5 + 15 + 12 + 2$). What's most fascinating about the Durer square is the bottom row: 1514 was the year the painting was made, and 4 and 1 correspond to the letters D and A, the initials of Albrecht Durer.

Today, we now know that there are exactly 7,040 fourth order normal magic squares and 2,202,441,792 fifth order normal magic squares in existence. What's most intriguing is that, like most other famous mathematicians of the time in our past, results were found without the modern-day computing and research power we have access to today. Algorithms were developed to create magic squares of \emph{any} size, without ever knowing how many there really are. This paper seeks to unveil that question: How to develop an algorithm that creates magic squares while systemically showing the methods it uses.

\section{Related Work}
For any normal magic square of order $n$, by careful algebra one can show that $$ \mu = \frac{n(n^2 + 1)}{2}$$ which, for the case of $n = 4$ is the number $34$. In a 1917 publication of \emph{Amusements in Mathematics}, Henry Dudeney, an English mathematician specializing in logic puzzles and mathematical games, used this important fact and found a way to classify all 7040 fourth order normal magic squares into 12 distinct configurations (coincidentally named Types I - XII). His classification was developed by considering the fact that the numbers from $1 .. 16$ contain 8 pairs of numbers that sum to $17 = \frac{\mu}{2}$, and diagrams can be drawn for each magic square to show these pairings (Fig 1.)
\begin{figure}[h]
\centering
\includegraphics[width=0.8\linewidth]{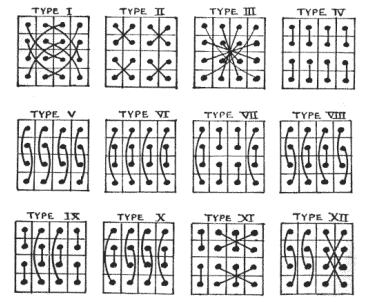}
\caption{The 12 Dudeney Classification Types}
%\label{Something else}
\end{figure}

These 12 classifications are independent of orientation, as many of the diagrams are not symmetric about some axis of reflection or rotation. Given this classification scheme, Dudeney later showed precisely how many of the 7040 magic squares were of each type
\begin{table}[h]
\scalebox{0.95}{
\begin{tabular}{c c c c c c}
\textbf{Type I}  & \textbf{Type II} & \textbf{Type III} & \textbf{Type IV} & \textbf{Type V} & \textbf{Type VI}   \\
384              & 384              & 384               & 768   & 768 & 2432             \\
& & & & &  \\
\textbf{Type VII} & \textbf{Type VIII} & \textbf{Type IX} & \textbf{Type X}  & \textbf{Type XI}  & \textbf{Type XII} \\
448               & 448 & 448              & 448              & 64                & 64                
\end{tabular}}
\end{table}

Curiously, one may note that Type I-III have an equivalent number of magic squares, as do Type VII-X and Type XI-XII. Looking back to Fig. 1, it also seems that these groups have similar structures in the configurations of half-sums. By 1948, in a paper \emph{Determinants of Fourth Order Magic Squares} published in the American Mathematical Monthly, C.W. Trigg proved that these similarities between certain Dudeney types existed by calculating their determinants. Trigg categorized the 12 Dudeney types into 4 distinct groups (unexcitingly named Type A-D). Indeed, Type A was the set of Type I-III, Type B the set of Type IV-VI, Type C the set of Type VII-X, and Type D the set of Type XI-XII. 

Furthermore, since the magic squares in these groups shared determinants, Trigg showed that certain row and column operations could be performed on magic squares preserving the determinant, but reordering the configuration of half-sums within that group. For instance, by interchanging the second and third rows, then interchanging the second and third columns of a Type II magic square, the result would be a Type I magic square (both within the Type A group). Alternatively, interchanging the third and fourth rows, then interchanging the third and fourth rows of a Type III magic square results in a Type I magic square as well.

\begin{figure}[h]
\centering
\includegraphics[width=0.8\linewidth]{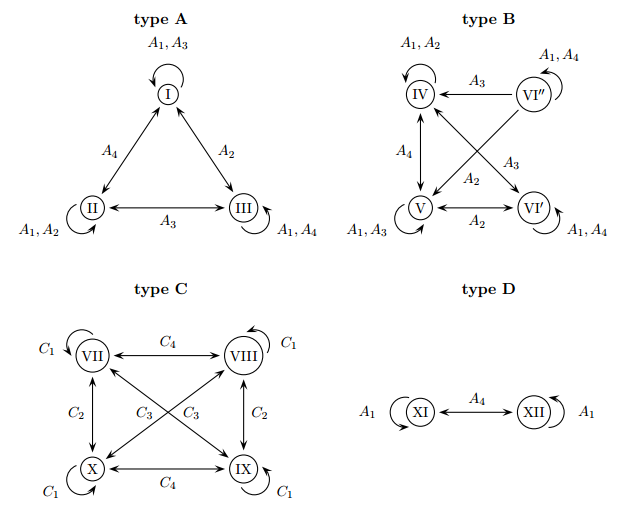}
\caption{Trigg Classification and Staab Transformations}
\end{figure}

Further research by Peter Staab in 2010 investigated the Trigg groups to find if additional operations existed within the groups to reconfigure one Dudeney type magic square to another, and if operations existed which did not affect the configuration of the square. Interestingly, he noted that Type VI needed to be split within the Trigg Type B into two distinct subsets: Type VI' and Type VI''. This was due to some of the Type VI magic squares having the additional property that a broken diagonal also summed to $\mu = 34$. After this reclassification, he found the following transformation groups existed for the Trigg classification (Fig 2.):
\begin{description}
\item[$A_1$] indentities and even bisymmetric transformations
\item[$A_2$] semisymmetric transformations
\item[$A_3$] broken-diagonal (translation) transformations
\item[$A_4$] odd bisymmetric and rotational transformations
\item[$C_1$] identity transformations
\item[$C_2$] odd bisymmetric transformations
\item[$C_3$] even bisymmetric transformations
\item[$C_4$] rotational transformations
\end{description}

Lastly, although tangential to the purpose of this paper, it is worth mentioning one key result from 1998. K. Pinn and C. Wieczerkowski approximated that, through the use of Parallel Tempering Monte Carlo simulations, the total number of sixth order normal magic squares is $(1.4196 \pm 0.00128) \cdot 10^{20}$. Implications from these findings suggest that modern day computing power may still not be fast enough to work in the space of sixth order magic squares.

\section{Problem Statement}
Given that much work has been done in categorizing the fourth order magic squares by their symmetries and pairwise similarities, the next logical step would be to investigate these symmetry groups to conclude which of the 7040 magic squares are necessary to \emph{generate} the fourth order space. We make the following definitions:
\newline \newline \textbf{Def:} \emph{An $n\times n$ \textbf{permutation matrix} is a square matrix whose entries are $\{0,1\}$ such that $1$ is in each column and each row of the matrix exactly once.}
\newline \newline \indent A permutation matrix $P$ acting on a matrix $A$ of similar dimension essentially reorders either the rows or columns of that matrix dependent on whether it is a left- or right-multiplication. For example, Consider the following matrix multiplication:
$$\left[ \begin{array}{ccc}
1 & 0 & 0\\
0 & 0 & 1\\
0 & 1 & 0
\end{array} \right] \left[ \begin{array}{ccc}
1 & 2 & 3 \\
4 & 5 & 6 \\
7 & 8 & 9
\end{array} \right] =\left[ \begin{array}{ccc}
1 & 2 & 3 \\
7 & 8 & 9 \\
4 & 5 & 6
\end{array} \right] $$
Note that the second and third rows of this matrix are swapped after being multiplied by the permutation matrix on the left. Thus $PA$ permutes rows and $AP$ permutes columns, while $PAP$ permutes both rows and columns.
\newline \newline \textbf{Def:} \emph{The \textbf{symmetry group} of a set of magic squares $G$ is the set of all pairs of permutation matrices $\{P_i,P_j\}$ in $\mathcal{M}^{n \times n}$ such that for any magic square $A$ in the set $G$, $P_i AP_j$ and $P_i A^{T}P_j$ is also a magic square and in the set $G$. Often times this will be called the \textbf{transformation group} of $G$.}
\newline \newline \textbf{Def:} \emph{Two magic squares $A$ and $B$ are \textbf{symmetric} if and only if $B$ can be written as $P_iAP_j$ or $P_iA^TP_j$ for a pair of permutation matrices $\{P_i,P_j\}$.}
\newline \newline \indent As an example, in \textbf{Fig. 2} the set $G$ can be thought of as Trigg's type A, and the transformation group of that set is any pair of permutation matrices from the set $A_1, A_2, A_3$, or $A_4$. Clearly any magic square in that set will be again be a magic square in the set if any of the former transformations are applied.
\newline \newline \textbf{Def:} \emph{The \textbf{order} of a transformation group $T$ of $G$ (denoted as $T(G)$) is the cardinality of the set, i.e., $|T(G)|$.}
\newline \newline \textbf{Def:} \emph{A magic square $A$ of a set of magic squares $G$ is a \textbf{generator} of $G$ if \begin{enumerate} \item The transformation group of $G$ applied to $A$ generates a subset of the elements of $G$ and \item For any two generators $A_i$, $A_j$ of $G$, the subsets of  $G$ generated by $A_i$ and $A_j$ are pairwise disjoint. \end{enumerate}} 
\bigskip The choice of the word generator in this case is not used in the traditional sense for algebra. A set of generators of a group of magic squares $G$ is simply the set of magic squares that are pairwise asymmetric and together generate the entire set $G$ using the transformation group of $G$. With these definitions in place, we are ready to approach the question to be answered:
\newline \newline \textbf{\emph{Of the fourth order normal magic squares, how many unique generators are necessary and sufficient to generate the entire space, and what are the associated transformation groups for these generators?}}

\bigskip As mentioned in the previous section, progress has been made towards the answer; however, the lower bound for generators itself is unconfirmed and the categorization of transformation groups is incomplete. Successfully determining a lower bound and completing the transformation groups would provide for two key results: \begin{enumerate} \item That there exists a set of magic squares for the n-th order that describe the entire space and \item That this space can be described by a categorical enumeration of transformation groups imposed on the space. \end{enumerate}
\indent Since the topology of matrices is well-defined and inherit many nice algebraic properties, the above results would hold for the space of magic squares of any order. This suggests that such a finding would be a methodical approach towards identifying the properties of higher order magic squares with relative ease.
\section{Methods}
$$ \begin{array}{|c|c|c|c|}
\hline 
a & b & c & d \\
\hline
e & f & g & h \\
\hline 
i & j & k & l \\
\hline
m & n & o & p \\
\hline
\end{array} $$
\subsection*{Theory}
As mentioned previously, a normal magic square $A$ of the fourth order has the following properties:
\begin{itemize}
\item $a_{i,j} \in A \implies a_{i,j} \in \{1, \ldots, 16\}$
\item Each $a_{i,j}$ distinct
\item $\mu = \frac{n(n^2 + 1)}{2} = 34$
\item Each row, column, and trace sum to $\mu$
\end{itemize}
The last two properties will be investigated first. Let's consider the generic $4 \times 4$ matrix from above that inherits these properties. We may notice that this produces a total of 10 linear constraints, namely
\begin{enumerate}
\item $a + b + c + d = 34$
\item $e + f + g + h = 34$
\item $i + j + k + l = 34$
\item $m + n + o + p = 34$
\item $a + e + i + m = 34$
\item $b + f + j + n = 34$
\item $c + g + k + o = 34$
\item $d + h + l + p = 34$
\item $a + f + k + p = 34$
\item $d + g + j + m = 34$.
\end{enumerate}
However, the solution to the system of simultaneous linear equations above finds that of the 16 variables in the matrix, only 7 are independent while the remaining 9 are necessarily defined by their dependence. That is to say, a fourth order normal magic square can be uniquely defined by only 7 entries, thus lowering the dimensionality of the problem greatly. For example, assume the entries $a, b, c, e, f, g, i$ are independent (the choices for these being independent will be examined later), then

\begin{description}
\item $d = 34 - a - b - c$
\item $h = 34 - e - f - g$
\item $j = 2a + b + c + e - g + i - 34$
\item $k = 68 - 2a - b - c - e - f - i$
\item $l = f + g - i$
\item $m = 34 - a - e - i$
\item $n = 68 - 2a - 2b - c - e - f + g - i$
\item $o = 2a + b + e + f - g + i - 34$
\item $p = a + b + c + e +i - 34$.
\end{description}

This does not necessarily suffice in showing that arbitrary numbers chosen for the entries $a,b,c,e,f,g,i$ will guarantee a normal magic square. It is certainly true that the constraints above will hold so that each row, column, and diagonal sum to 34, but often times (at least with an arbitrary arrangement), entries may be duplicated, and others could be negative, 0, or larger than 16. However, if we condition the above constraints using the first two properties from above, we yield an important restriction on the arrangement of numbers for the independent cells: They must be configured in such a way as to obey the laws that every cell must contain a distinct number from the set $\{1, \ldots, 16\}$.

\subsection*{Implementation}
Determining the arrangements which do not follow this law requires an algorithmic approach, due to having no prior knowledge of which arrangements yield a normal magic square. Because we are tasked with many constraints, of which require very specific parameters (each entry being distinct), a backtracking algorithm that used constraint propagation seemed to be the most logical implementation, considering the analysis was focused on creating an efficient, dynamic system. Furthermore, backtracking has the property of quickly withdrawing from partial solutions which are guaranteed to not resolve as a full solution---for instance, when one of the constraints won't be met since the number necessary to be assigned has already previously been assigned. These assignments are carried out by the constraint propagation, which accounts for all properties listed in the previous section.

The algorithm begins by placing arbitrary numbers along the top row into slots $a, b, c$. The numbers in slots $a,b$ are considered "free" since they do not propagate any constraints, however, after $a$ and $b$ have been assigned, $c$ must be assigned in such a way that: $34 - a + b + c$ is a valid number (1 to 16) that has not been propagated or previously assigned. For instance, if $a = 1, b = 3$ then the assignment $c=15$ would be invalid, since $34 - a - b - c = 15$ which is previously assigned to $c$.

\begin{figure}[!h]
\centering
\includegraphics[width=0.8\linewidth]{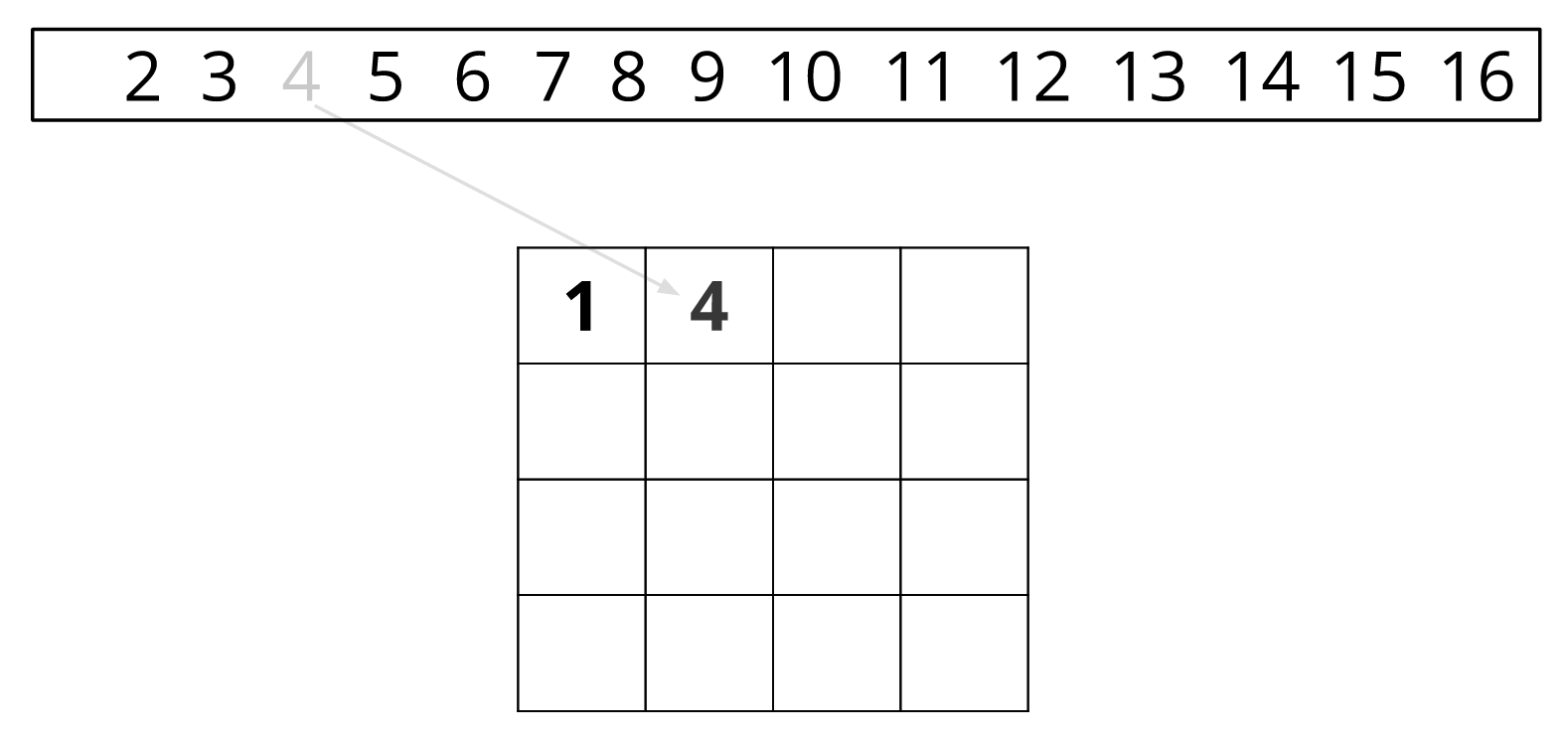}
\caption{Backtracking Algorithm by placing numbers into matrix}
\end{figure}

Should a valid assignment occur for $a,b,c$, then the number resulting from $34-a-b-c$ should be removed from the list of available terms for the backtracking algorithm to test. This is necessary, as without it, deeper iterations would test this number as an assignment despite it being necessarily assigned prior to $d$. After computing the first row, the algorithm computes the first column, having already registered the assignment for $a$. The choice for this is not so trivial: The cells $e, f$ are also considered "free" since they will not propagate any constraints. Because of this, "free" cells check for every valid assignment from the list of available numbers. The goal would then be to avoid "free" cells until absolutely necessary. Assigning the first column requires only one "free" cell, which in turn requires much less complexity.

\begin{figure}[!h]
\centering
\includegraphics[width=0.8\linewidth]{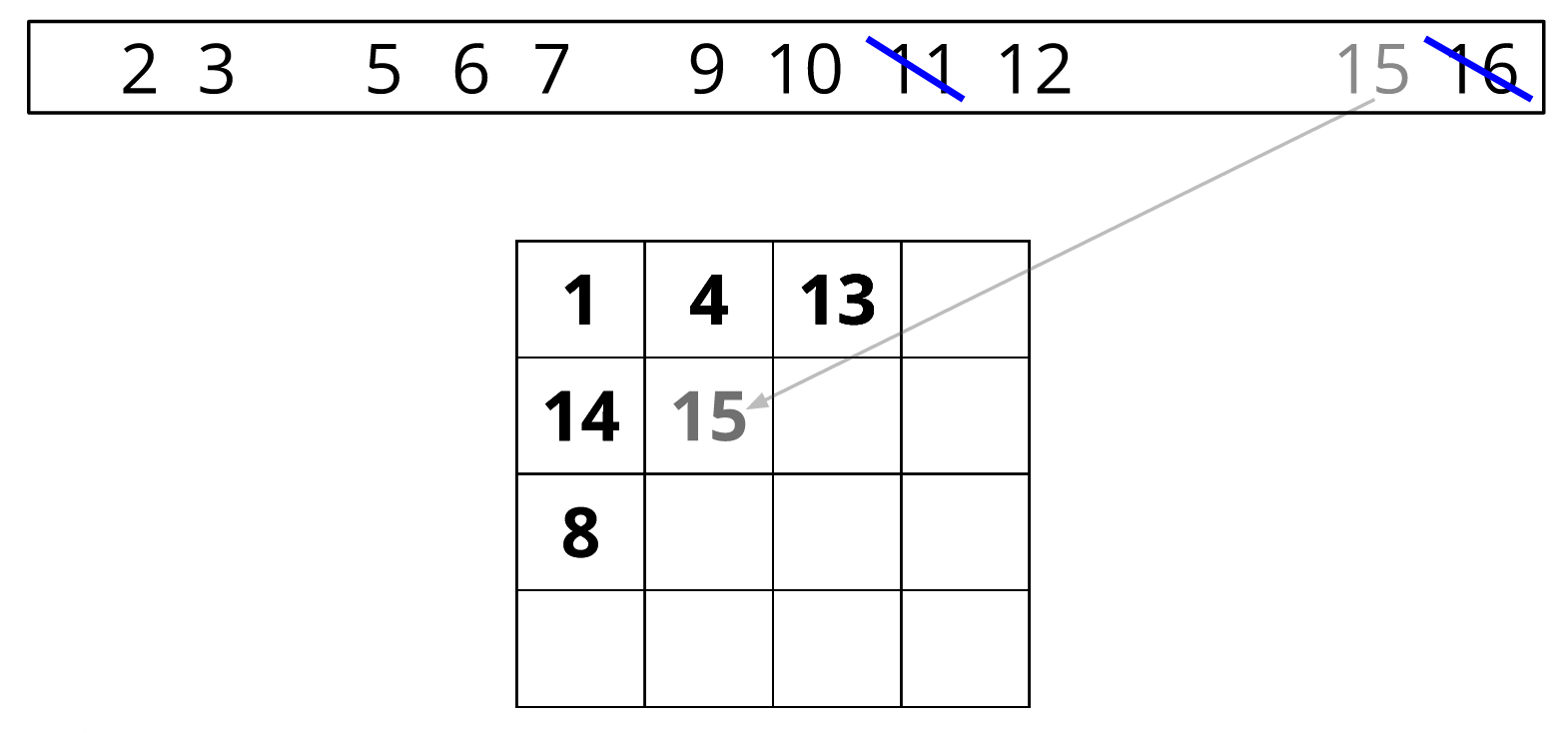}
\caption{Constraint propagation of values resulting in $\mu = 34$}
\end{figure}

After completing the first column, the last and final step of the backtracking algorithm is to compute the assignment for $f, g$, as they are the last two remaining cells necessary to be assigned to thusly define the square unique to these 7 numbers. From above, $f$ is a free cell, so values are arbitrarily assigned from the available list. The cell $g$ requires the most effort: Since this is the final cell, the number assigned here must be such that all constraints that contain $g$ are satisfied, and the dependence values are available. If either of these two assumptions fail, then $g$ cannot be assigned and the backtracking algorithm withdraws from this arrangement of numbers. 
\begin{figure}[!h]
\centering
\includegraphics[width=0.8\linewidth]{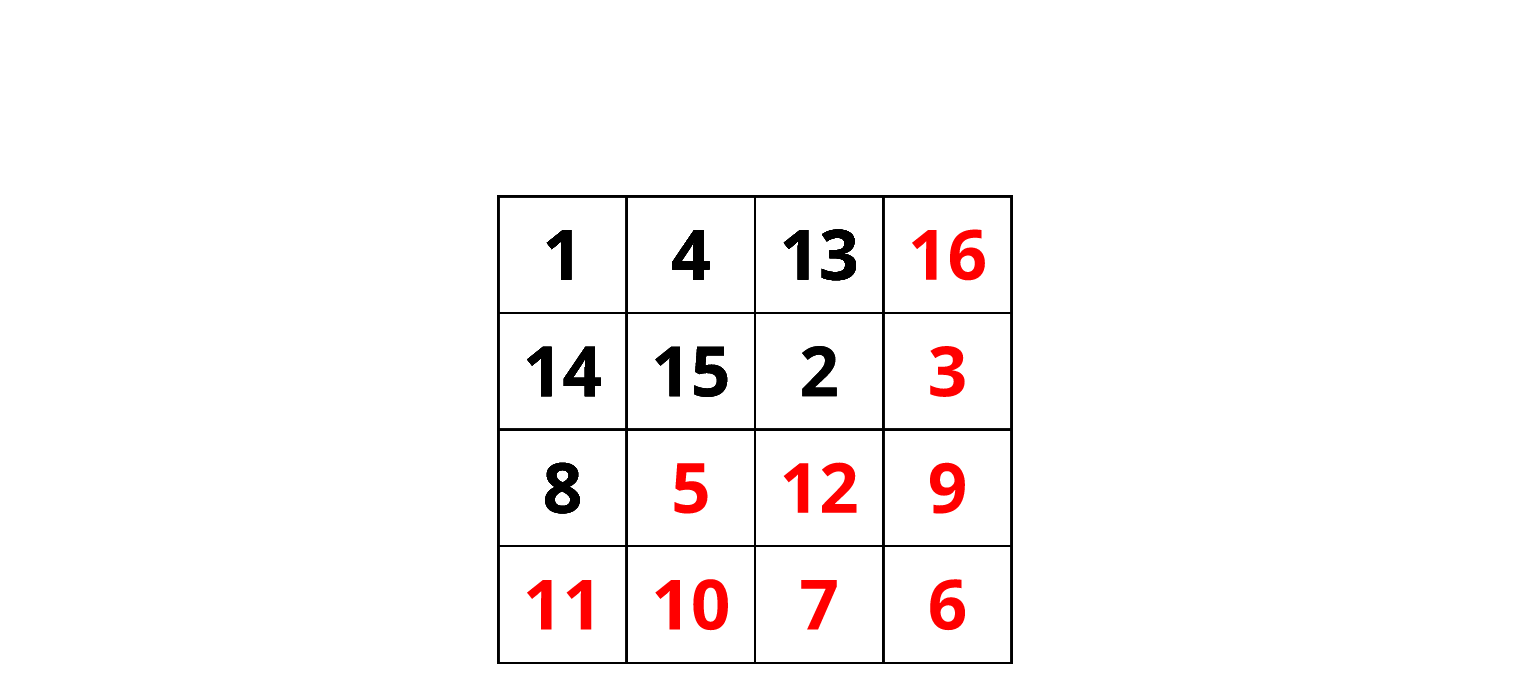}
\caption{Resultant magic square given by 7-dimension basis from above}
\end{figure}

Once the algorithm successfully generates a normal magic square using the independent 7 numbers for the cells (Fig. 5), a simple supervised machine learning technique is employed to find the pairs of 17 (Dudeney's classification scheme) in the magic square that can uniquely classify it. This is done by creating the classification scheme prior to using the diagram in Fig. 1. 
\begin{figure}[!h]
\centering
\includegraphics[width=0.8\linewidth]{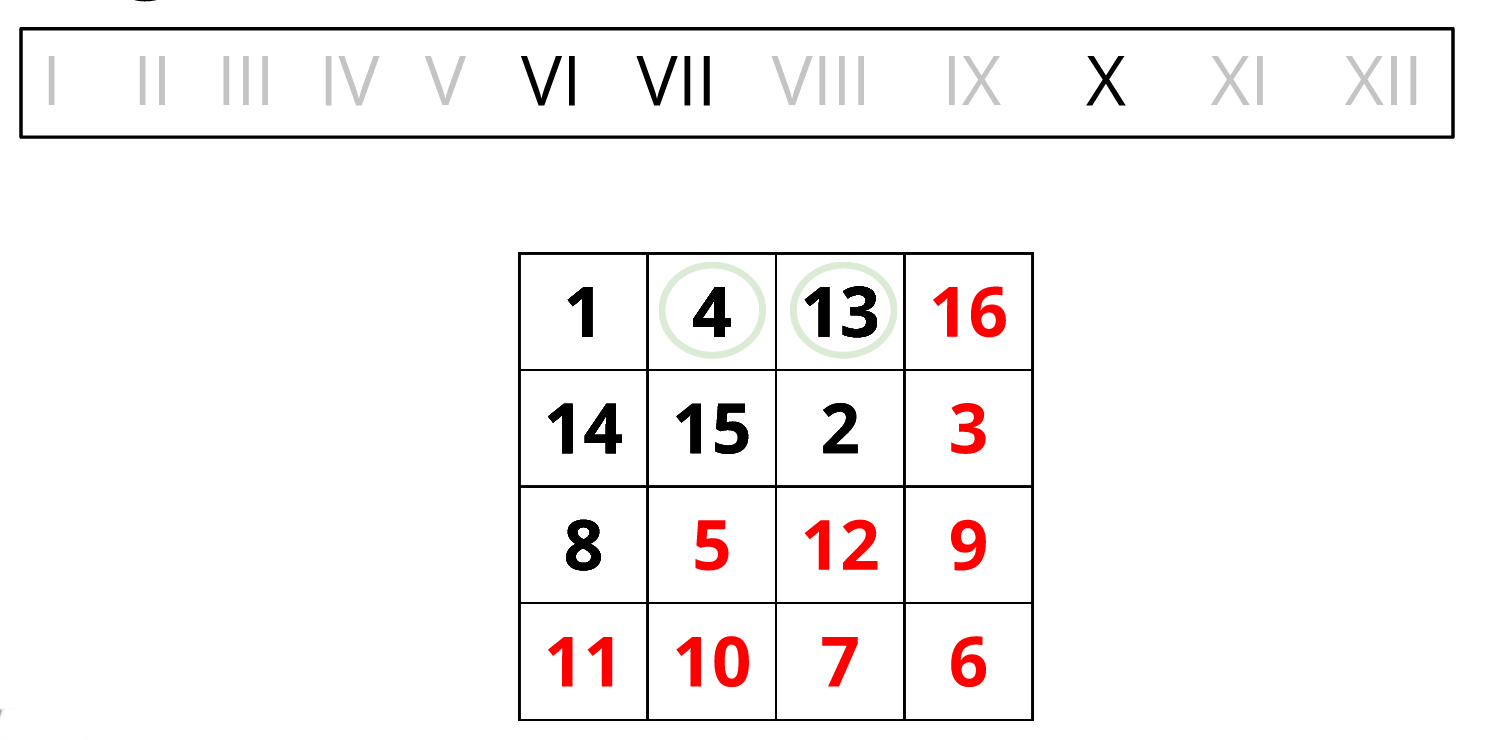}
\caption{Categorization by finding sums of $\mu /2 = 17$}
\end{figure}

Then, using the 7 numbers that define the magic square, the program naively iterates across the set to find the pairs. This is surprisingly efficient, since it does not have to iterate across all 16 numbers to find pairs, and the classification scheme guarantees that only two pairs are necessary to determine the Dudeney Type.

\begin{figure}[!h]
\centering
\includegraphics[width=0.8\linewidth]{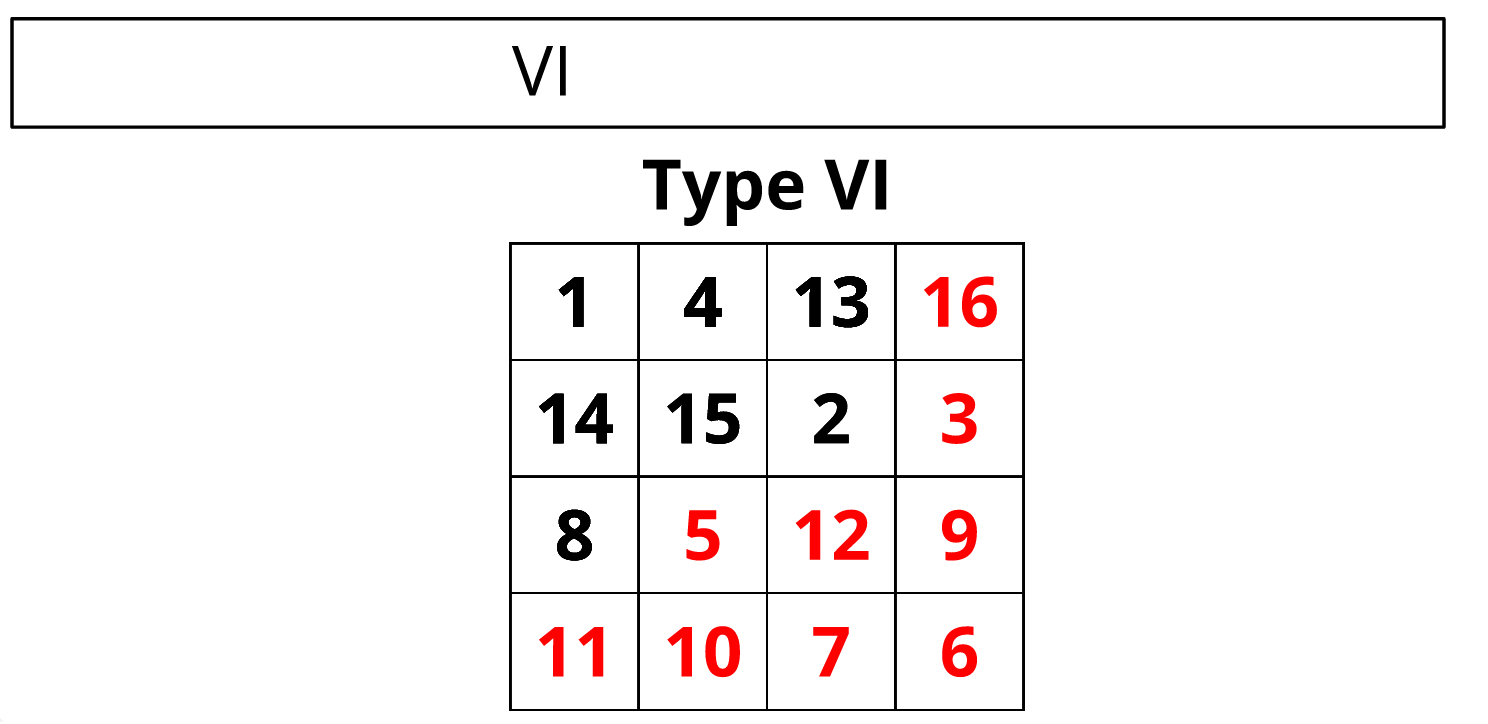}
\caption{Resultant classification type given by configurations of sums}
\end{figure}

After each square is generated and classified, the next step was to classify each square, as an extension, by their transformation groups. The properties necessary for this included their innate symmetries (such as how many broken diagonals also added to $\mu$), which Trigg group they were in, and also by using previous research of the permutation sets Staab had proven.

\begin{figure}[!h]
\centering
\includegraphics[width=0.8\linewidth]{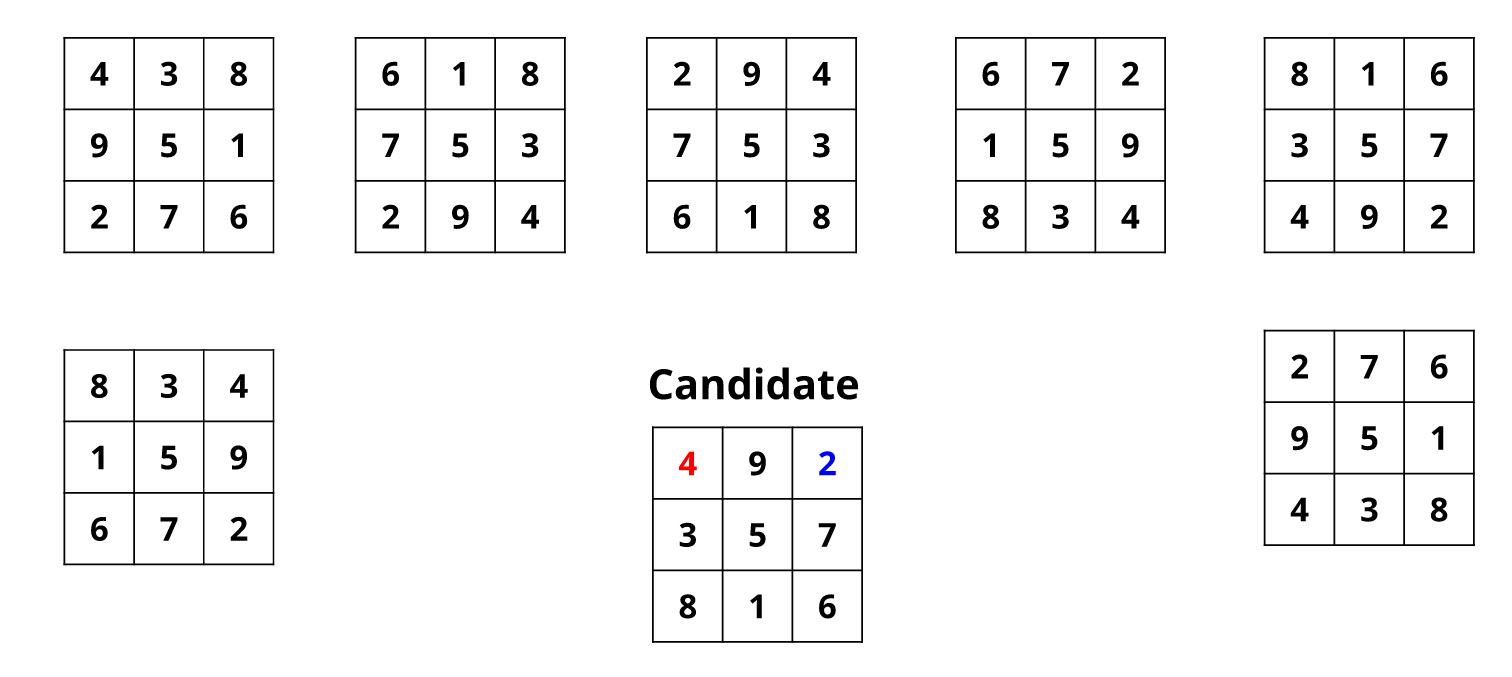}
\caption{Candidate magic square from symmetric transformation group}
\end{figure}
\vfill\eject
Finally, an arbitrary magic square was chosen from each transformation group within the Trigg type, and had applied to it every set of permutation matrices defined on that magic square. As the permutations were applied, magic squares that were symmetric were deleted from the group, until the transformation group had been exhausted. If more than one magic square remained in the group, the process was repeated with another magic square from the group (which clearly defined another subset of magic squares which it generates) until the minimum was reached.
\begin{figure}[!h]
\centering
\includegraphics[width=0.8\linewidth]{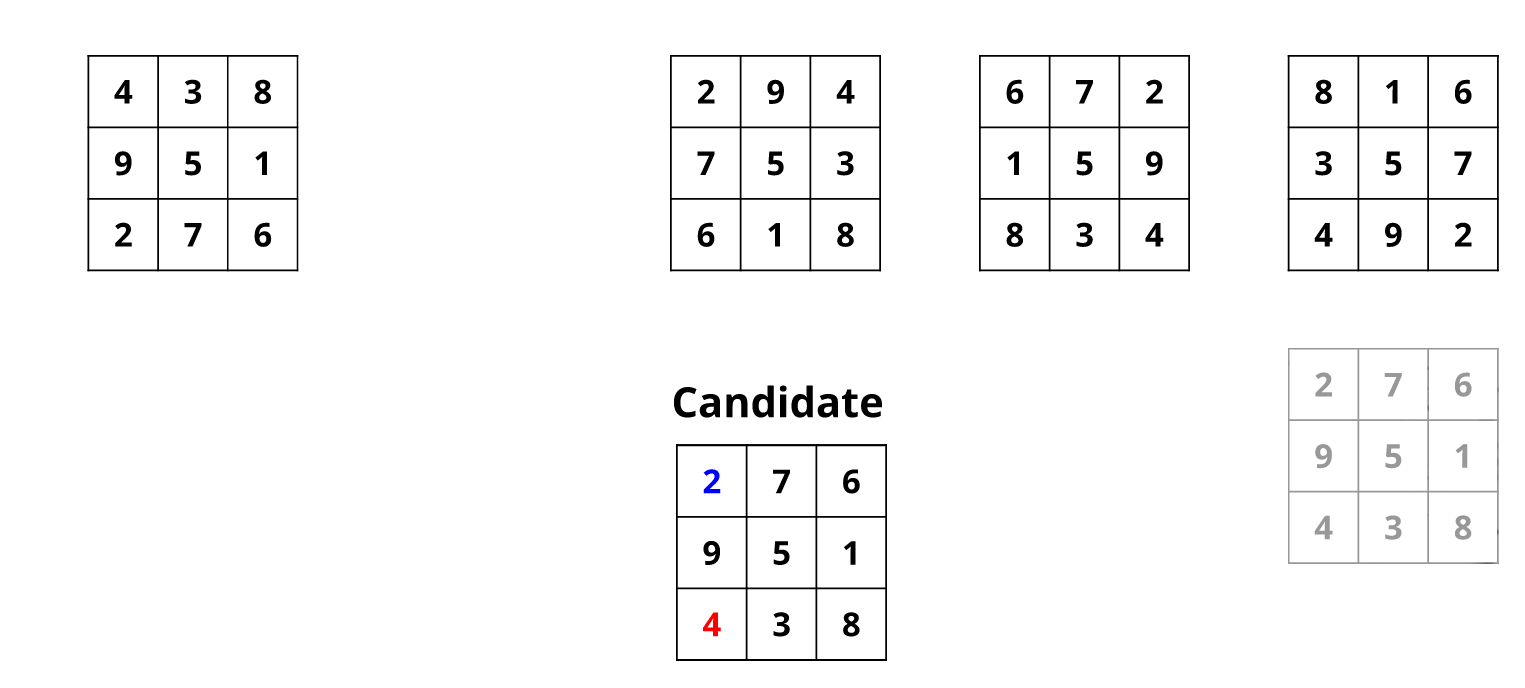}
\caption{Removal of magic squares in the group based on rotations and transposition}
\end{figure}

An example of this transformation group search can be seen for the third-order case in Fig. 8, 9, and 10. If we consider the dihedral group (rotations, reflections, transpositions) as the transformation group for this set, then you can see that as the Candidate square rotates about, transposes, and reflects about each axis, it becomes symmetric with every other magic square in the set, thus deleting them from the group. Note that since there are only 8 normal magic squares of order 3, that the candidate magic square is indeed the sole generator of the entire group.
\begin{figure}[!h]
\centering
\includegraphics[width=0.8\linewidth]{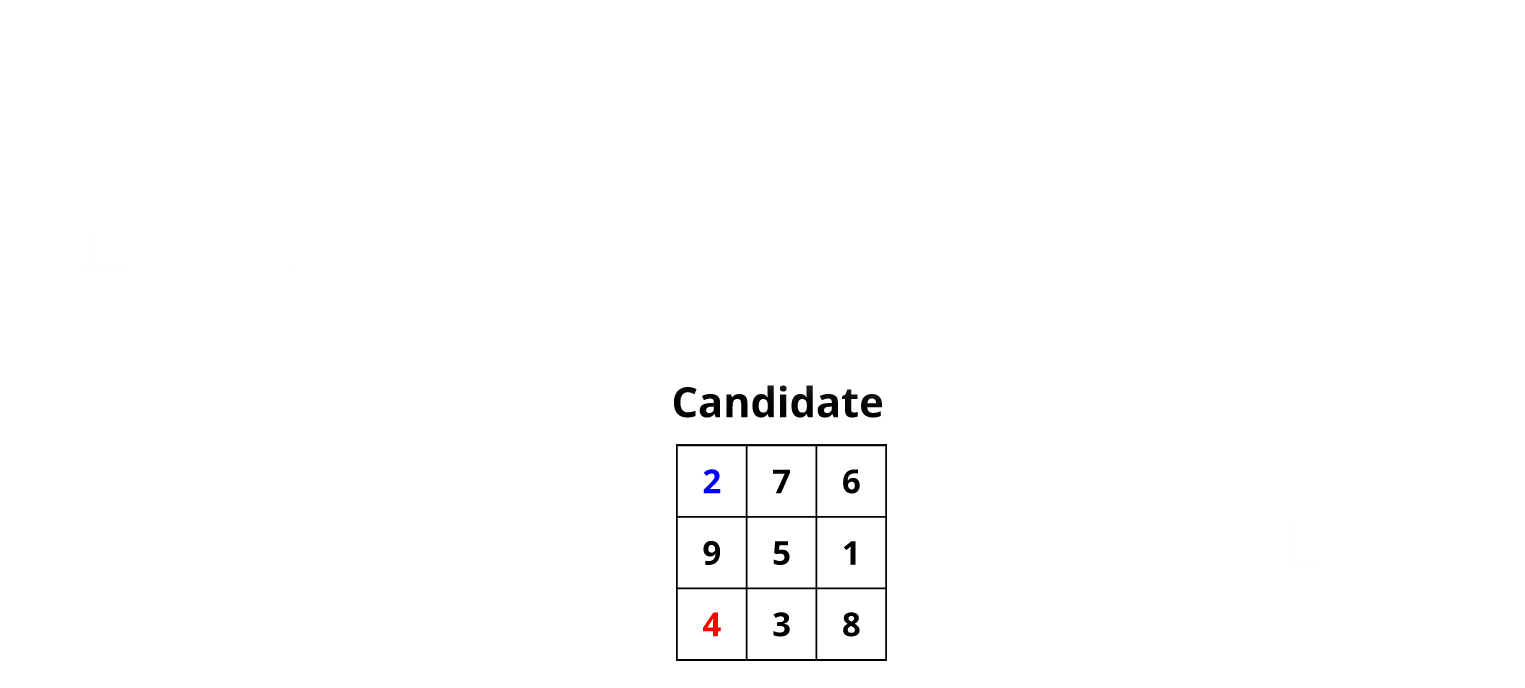}
\caption{Resultant unique generator for order 3 magic squares}
\end{figure}

\section{Results}
Even though there are 7040 magic squares, results found that only 95 magic squares are necessary and sufficient to generate the entire space of fourth order. What was immediately perplexing about this number was how it was distributed among the different Trigg Types. For the Type B and Type C magic squares, the resultant number of generators for the group did not evenly divide the enumeration of the magic squares classified as that type; the expectation was that a transformation group would provide a reduction of the set by some factor, but evidence suggested that additional information was necessary. 

\subsection*{Trigg Type A}
The magic squares that exist in this group, Dudeney Types I-III, have the additional property that each broken diagonal also sums to $\mu = 34$. Because of this property, many symmetries are introduced into this group. Analysis of the results found that, from the 1152 magic squares classified as Type A, only 3 generators are necessary to generate the entire group (Fig. 11).
\begin{figure}[h]
\centering
\includegraphics[width=0.8\linewidth]{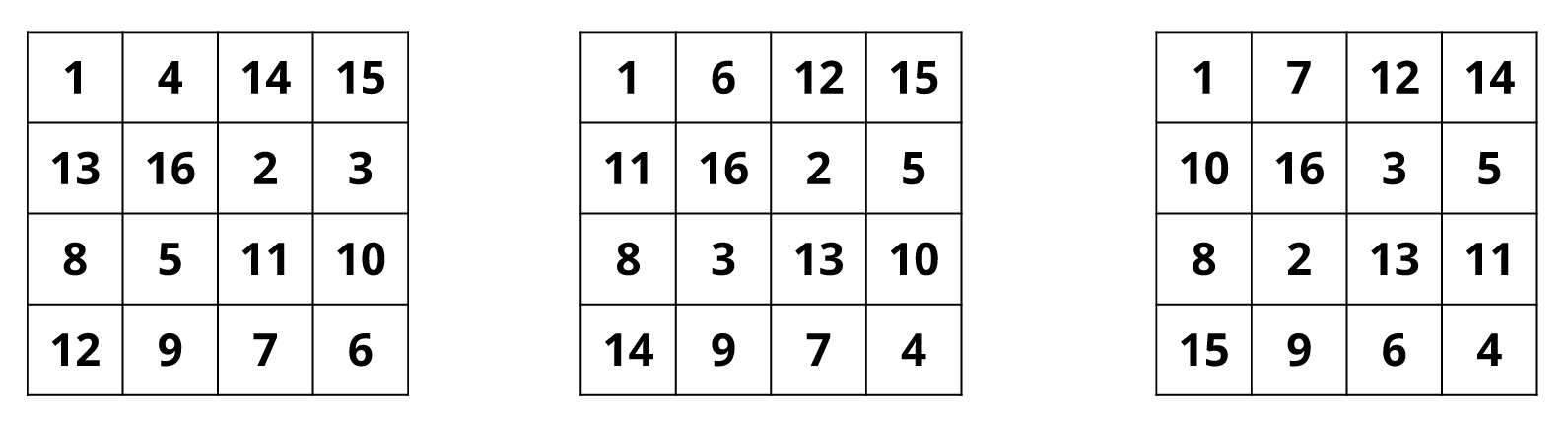}
\caption{Generators for Type A}
\end{figure}

Thus each generator from this set has a transformation group of order 384. Furthermore, we can claim that the entire Type A space is uniquely defined by these 3 generators and the transformation group of order 384. One observation worth mentioning is how similar each of the generators are; Only 8 of 16 cells differ between each generator, and yet there is no arrangement of permutations in the confines of the space to transform one to the other.

\subsection*{Trigg Type B}
Analysis of the Dudeney Types IV-VI provided the most interesting results, primarily in the fact that so little was previously understood. It was mentioned earlier that Type VI had the additional property that one of its broken diagonals also summed to $\mu = 34$. The assumption, then, would be that this group would undoubtedly have generators belonging to two separate transformation groups.
\begin{figure}[h]
\centering
\includegraphics[width=0.95\linewidth]{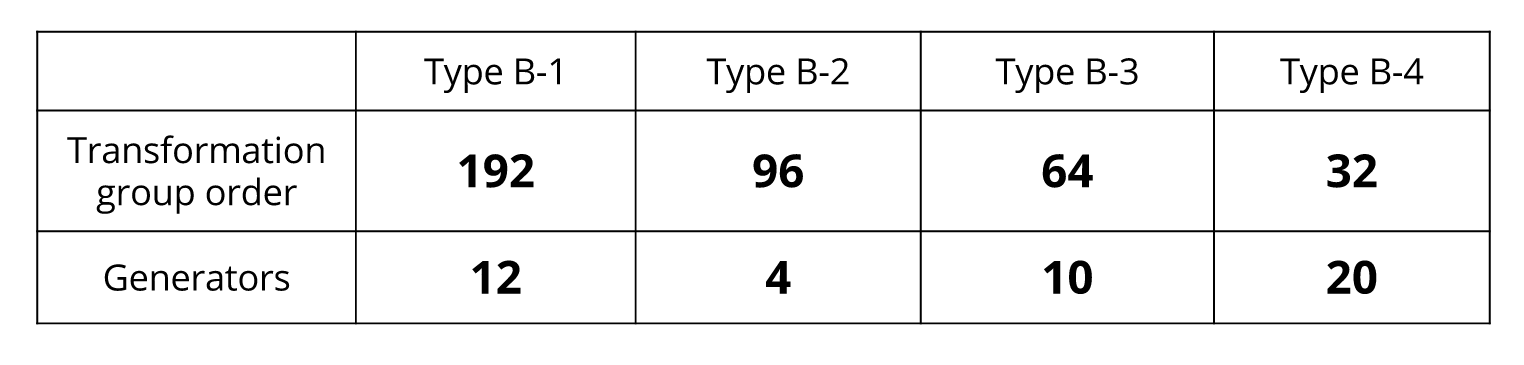}
\caption{Generators for Type B}
\end{figure}

This was found to be not entirely true (Fig. 12). In fact, hidden inside the group of 3968 magic squares classified as Type B was a total of 46 generators, each belonging to one of 4 different transformation groups. The transformation groups were found to have orders 192, 96, 64, and 32, with generators 12, 4, 10, and 20 respectively. Unlike Type A, Type B space could not be uniquely defined by its generators and a transformation group; however, if one were to further reduce Type B into 4 pairwise disjoint subsets (aptly named B-1, B-2, B-3, and B-4) containing only the generators of a particular transformation group, then we can indeed uniquely define each \emph{subgroup} by its generators and the transformation group. 

\subsection*{Trigg Type C}
While analysis of the Trigg Type B group provided intriguing results, Trigg Type C---Dudeney Types VII-X---had perhaps the most surprising result of all four types. Given that there were 1792 magic squares in Type C evenly distributed among all four Dudeney types, the assumption would be that the behavior of this group was much like Type A. 
\begin{figure}[h]
\centering
\includegraphics[width=0.6\linewidth]{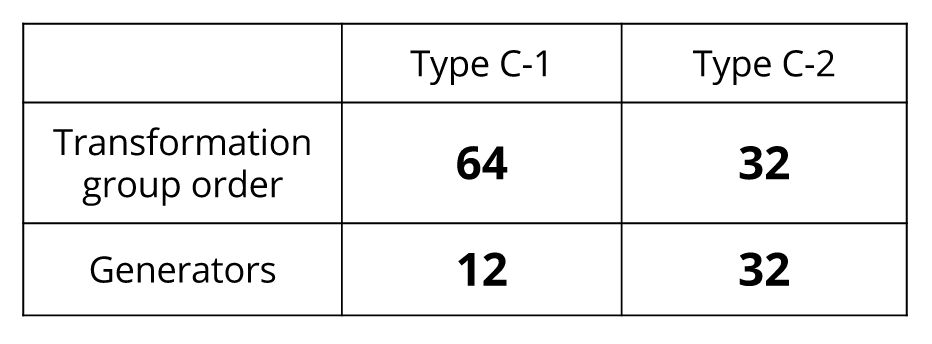}
\caption{Generators for Type C}
\end{figure}
In fact, it was found that there were 44 generators for this group, each belonging to one of 2 different transformation groups. The transformation groups were found to have orders 64 and 32, with generators 12 and 32, respectively (Fig. 13). Again, because of this result, the justification that Type C could be uniquely defined by a set of generators and a transformation group was invalid. By imitating the revision for Type B and splitting the Type C group into two disjoint subsets (C-1, C-2) categorized by separate generators and their transformation sets, we again arrive at a distinction of uniquely defining a group of magic squares based on those generators.

\subsection*{Trigg Type D}
With the small amount of magic squares classified Dudeney Type XI-XII in this group (128), the results were expectantly underwhelming and categorically likened to Type A. The exception in this case is that Trigg Type D magic squares lack many symmetries which would provide for a high-ordered transformation group (Fig. 1). 
\begin{figure}[h]
\centering
\includegraphics[width=0.8\linewidth]{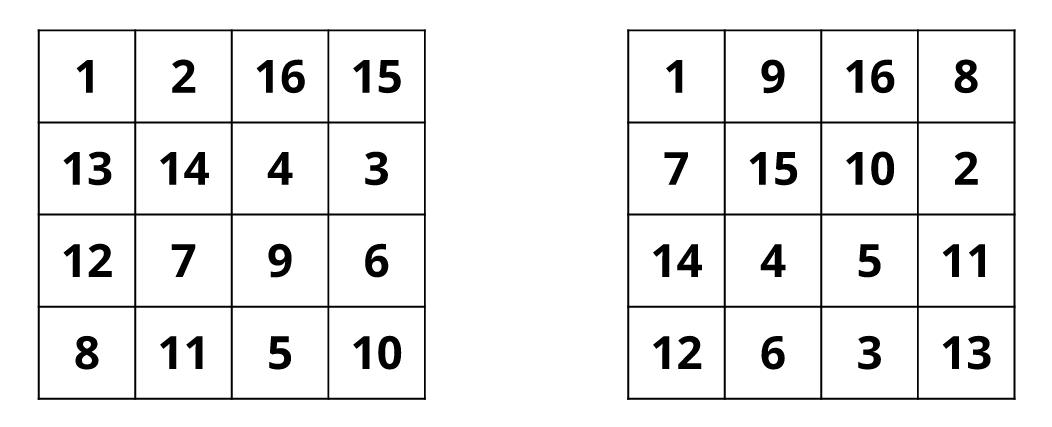}
\caption{Generators for Type D}
\end{figure}

It was found that only 2 generators were necessary to generate the space of Type D magic squares, and each generator had a transformation group of order 64 (Fig. 14). Much like Type A, without the need to reclassify the group, we can uniquely define the Type D magic squares by these two generators and the corresponding transformation group of order 64. Another observation worth mentioning related to the similarity of these generators: There's practically none, and yet due to the lack of symmetry within the group, no arrangement of permutations in the confines of this space generate one from the other.
\section{Discussion}

The results above motivate a theoretical discussion on the potential of the methodology scaled to higher orders. Firstly, it was found that very minimal information was necessary to completely describe the entire space of fourth order magic squares. From the previous sections, it was noted that the problem in itself could be lower-dimensionalized to a basis of only 7 out of the 16 numbers that define a magic square. Further, only 95 generators are sufficient to generate the 7040 magic squares for order four. This implies, that given the knowledge of Trigg classification and transformation groups, less than 700 bytes of information ($7*95$ and assuming each number requires one byte of storage) is required. 

Since we know \emph{some} properties about higher order magic squares, such as transformation groups of odd order have the same behavior as transformation groups of one less order (even), we can hypothesize that this scheme is adaptable--at least for the fifth order case. This is further emphasized by the fact that the methods used in this analysis were efficient and easy to implement. With careful linear algebra, one can show that the fifth order magic squares can be lower-dimensionalized to a basis of only 13 out of 25 numbers that define that magic square. The implications from this suggest that a backtracking algorithm with constraint propagation may still be suitable and efficient enough with modern computing power to produce meaningful results.

Lastly, there is consideration to be had to appropriate a better approximation (or more concrete bound) on the total number of sixth order magic squares. While a previously mentioned result from 1998 showed an approximation for this number, certain applications from this analysis could perhaps be abstracted to find the theoretical transformation groups of the set, and determine from those what the generators of those sets would look like. While there is no information currently in literature suggesting any formulation of "Trigg Types" for these magic squares, one could again theorize their configurations based on small simulations of this method and applying transformation groups to this small subset.
\vfill\eject
\section{Conclusion}
The finding of 95 generators for the entire set of 7040 fourth-order normal magic squares improves the previous lower bound of 220 found in literature[2]. Furthermore, while many of the transformation groups were correctly categorized by Staab [5], results show that previously unknown transformation groups exist in the space, and that the Trigg classification requires revisions to reflect the true underlying structure. With these two results combined, by virtue of the methods used in this analysis, such a construction can be developed with ease for fifth order magic squares which has previously been limited by inefficient methods. With this paper as reference, a curious individual in the field of recreational mathematics may use these tools to yield many fascinating results about fifth and sixth order normal magic squares. 
%\vfill\eject


%\section*{Acknowledgment}




%\bibliographystyle{IEEEtran}

\begin{thebibliography}{1}

\bibitem{IEEEhowto:dudeney}
H. E. Dudeney, \emph{Amusements in Mathematics}, \hskip 1em London: Thomas Nelson, 1917

\bibitem{IEEE:howto:gardner}
M. Gardner, \emph{Mathematical Games} Scientific American iss. 234, pp. 118-122, 1976

\bibitem{IEEEhowto:pinn}
K. Pinn, C. Wieczerkowski, \emph{Number of Magic Squares From Parallel Tempering Monte Carlo} \hskip 1em arXiv.org, 1998

\bibitem{IEEEhowto:schroeppel}
R. Schroeppel, \emph{Appendix 5: The Order 5 Magic Squares} pp. 1-16, 1973

\bibitem{IEEEhowto:staab}
P. Staab, \emph{The Magic of Permutation Matrices: Categorizing, Counting, and Eigenspectra of Magic Squares}\hskip 1em \relax arXiv.org, 2010

\bibitem{IEEEhowto:trigg}
C. W. Trigg, \emph{Determinants of Fourth Order Magic Squares} \hskip 1em \newline American Mathematical Monthly iss. 55, pp. 558-561, 1948

\end{thebibliography}

\end{document}